\documentclass{article}
\usepackage{sw20bams}
\usepackage{latexsym}
\usepackage{amsmath}
\usepackage{graphicx}
\usepackage{amsfonts}
\usepackage{amssymb}
\newtheorem{theorem}{Theorem}
\newtheorem{acknowledgement}[theorem]{Acknowledgement}

\newtheorem{corollary}[theorem]{Corollary}

\newtheorem{definition}[theorem]{Definition}

\newtheorem{lemma}[theorem]{Lemma}

\newtheorem{problem}[theorem]{Problem}
\newtheorem{proposition}[theorem]{Proposition}
\newtheorem{remark}[theorem]{Remark}

\begin{document}

\author{Andrey Todorov\\University of California,\\Department of Mathematics,\\Santa Cruz, CA 95064.\\Bulgarian Academy of Sciences\\Institute of Mathematics\\Sofia, Bulgaria}
\title{Ray Singer Analytic Torsion of CY Manifolds II.}
\maketitle
\begin{abstract}
In this paper we construct the analogue of Dedekind eta function for odd
dimensional CY manifolds. We \ use the theory of determinant line bundles. We
constructed a canonical holomorphic section $\eta^{N}$ of some power of the
determinant line bundle on the moduli space of odd dimensional CY manifolds.

According to Viehweg the moduli space of moduli space of polarized odd
dimensional CY manifolds $\mathcal{M}(M)$ is quasi projective. According to a
Theorem due to Hironaka we can find a projective smooth variety $\overline
{\mathcal{M}(M)}$ such that $\overline{\mathcal{M}(M)}\backslash$
$\mathcal{M}(M)=\mathcal{D}_{\infty}$ is a divisor of normal crossings. We
also showed by using Mumford's theory of metrics with logarithmic growths that
the determinant line bundle can be canonically prolonged to $\overline
{\mathcal{M}(M)}$ $.$ We also showed that there exists section $\eta$ of some
power of the determinant line bundle which vanishes on $\mathcal{D}_{\infty}$
and has a Quillen norm the Ray Singer Analytic Torsion$.$ This section is that
analogue of the Dedekind eta Function $\eta.$
\end{abstract}
\tableofcontents

\section{Introduction.}

\subsection{Analytic Torsion of Elliptic Curves and Kronocker's Limit Formula}

In case of elliptic curves the exponential of Ray Singer analytic torsion is
the Quillen norm of a non vanishing section of the determinant line bundle.
Kronecker limit formula states that the Ray Singer analytic torsion is exactly
the Quillen norm of the Dedekind eta function $\eta.$ One of the versions of
the proof of the Kronecker limit formula is based on the facts that log of the
regularized determinant of the flat metric on the elliptic curve is the
potential for the Poincare metrics and this determinant is bounded as a
function on the moduli space. For this proof of the Kronecker limit formula
see \cite{JT98}.

Kronecker limit formula can be interpreted as the existence of a section
$\eta^{24}$ of some power of the determinant line bundle whose Quillen norm is
exactly the 24$^{th}$ power of the Ray Singer Analytic torsion.

\subsection{Formulation of the Problem Discussed in the Paper}

\begin{problem}
\label{PI}Does there exists a section $\eta^{N}$ of some power of the
determinant line bundle over the moduli space of odd dimensional CY manifolds
whose Quillen norm is the N$^{th}$ power of the Ray Singer Analytic torsion?
\end{problem}

The positive answer to Problem \ref{PI} can be considered as a\ generalization
of the above mentioned relations between Ray Singer Analytic torsion and the
Quillen norm on the determinant line bundle from elliptic curves to odd
dimensional Calabi-Yau manifolds. Moreprecisely the positive answer of the
Problem \ref{PI} will be a generalizion of Kronecker's Limit Formula to higher dimensions.

\subsection{Description of the Ideas Used in the Paper}

The idea on which this paper is based, is very simple. The Quillen metric is
related to the spectral properties of the Laplacians acting on (0,q) forms in
case of K\"{a}hler manifolds. The main question is when the Ray Singer
analytic torsion is the Quillen metric of some holomorphic section of the
determinant line bundle.

One of the main results of the paper is the construction of a canonical
section of the determinant line bundle over the moduli space of odd
dimensional CY manifolds up to a constant whose absolute value is one. First
we prove that the determinant line bundle $\mathcal{L}$ as $C^{\infty}$ bundle
is trivial. The proof that the determinant line bundle $\mathcal{L}$ is a
trivial $C^{\infty}$ is based on the facts that for odd dimensional CY
manifolds we have dim $\ker($ $\overline{\partial}$ $\overline{\partial}%
^{\ast})=$dim $\ker($ $\overline{\partial}^{\ast}\overline{\partial})=1$ over
the moduli space and the Ray Singer Analytic Torsion I(M) is not a constant
and strictly positive function on the moduli space. The definition of the
Quillen metric on the determinant line bunlde $\mathcal{L}$ implies that
$c_{1}(\mathcal{L})=d(\partial(\log(I(M)).$ Therefore $c_{1}(\mathcal{L})$ is
an exact two form on $\mathcal{M}(M)$.

It is easy to see that the determinant line bundle is isomorphic to $\pi
_{\ast}(\Omega_{\mathcal{X}/\mathcal{M}(M)}^{2n+1}).$ This implies that we can
construct a non zero section of the determinant line bundle whose $L^{2}$ norm
pointwise is equal to one. Using that section we can construct a canonical
C$^{\infty}$ non vanishing section of the determinant line bundle up to a
constant whose Quillen norm is exactly the analytic Ray Singer torsion. We
will call such section det($\overline{\partial}).$

It is not difficult to see that the analytic torsion for even dimensional CY
manifolds is equal to zero and for odd dimensional CY manifolds it is
different from zero. On the other hand the index of the $\overline
{\partial\text{ }}$ operator on the complex of $(0,q)$ forms is zero for odd
dimensional CY manifolds therefore there exists a canonical non vanishing
section det($\overline{\partial})$ whose Quillen norm is exactly the Ray
Singer Analytic Torsion.

We will study the zero set of the canonical section $\det(\overline{\partial
})$ in this paper on some compactification of $\overline{\mathcal{M}^{\prime
}\mathcal{(}M)}$ such that $\overline{\mathcal{M}^{\prime}\mathcal{(}M)}$
$\backslash$%
$\mathcal{M}^{\prime}\mathcal{(}M)$ is a divisor of normal crossings. Viehweg
proved in \cite{W}\ that $\mathcal{M}^{\prime}\mathcal{(}M)$ is a
quasi-projective variety. Moreover it is a well known fact that the moduli
space of CY manifolds is obtained by factoring the Teichm\"{u}ller space by
subgroup of finite index in the mapping class group that preserves some
polarization on the CY manifolds, which according to Sullivan is an arithmetic
group. From the fact that the mapping class group is an arithmetic one can
find a subgroup of finite index in the mapping class group such that the
quotient of the Teichm\"{u}ller space by this group is a non-singular variety
$\mathcal{M(}M)$ which is a finite covering of $\mathcal{M}^{\prime
}\mathcal{(}M).$

Let $\overline{\mathcal{M(}M)}$ be some compactification of $\mathcal{M(}M)$
such that $\overline{\mathcal{M(}M)}\backslash\mathcal{M(}M)\frak{=D}_{\infty
}$ is a divisor with normal crossings. The divisor $\frak{D}_{\infty}$ will be
called the discriminant locus. It is easy to see that the determinant line
bundle $\mathcal{L}$ for odd dimensional CY manifolds is isomorphic to the
dual of the holomorphic line bundle $R^{0}\pi_{\ast}\Omega^{2n+1}.$ On the
line bundle $R^{0}\pi_{\ast}\Omega^{2n+1}$ we have a natural metric
$\parallel$ $\parallel^{2}.$ One can show that the metric $\parallel$
$\parallel^{2}$ has a logarithmic growth in the sense of Mumford. See
\cite{Mu}. From here we deduced that determinant line bundle $\mathcal{L}$ can
be prolonged in an unique way to a line bundle $\overline{\mathcal{L}}$ over
$\overline{\mathcal{M(}M)}$ using the metric $\parallel$ $\parallel^{2}.$

Since we proved in \cite{T1} that the Quillen norm of the canonical section
$\det(\overline{\partial})$ is bounded, it is not difficult see that we can
continue the C$^{\infty}$ section $\det(\overline{\partial})$ to a section
$\overline{\det(\overline{\partial})}$ of the line bundle $\overline
{\mathcal{L}}$ over $\overline{\mathcal{M(}M)}$ and whose zero set is
supported by $\mathcal{D}_{\infty}.$

\subsection{Formulation of the Main Result}

The main theorem of the paper is the following one:

\textbf{THEOREM.} There exists a holomorphic section $\overline{\eta}^{N}$ of
the line bundle $\left(  \overline{\mathcal{L}}\right)  ^{\otimes N}$ over
$\overline{\mathcal{M(}M)}$ whose zero set is supported by $\mathcal{D}%
_{\infty}$ and the Quillen norm of $\overline{\eta}^{N}|_{\mathcal{M}%
(M)}=I(M)^{N}.$

\subsection{Outline of the Proof of the Main Result}

The ideas of the proof are the following:

\textbf{Step 1. }We noticed in \cite{T1} that the analytic torsion is equal to
the determinant of the Laplacians det($\triangle_{\tau}$) of Calabi-Yau
metrics acting on functions, whose imaginary part has a fixed cohomology
class, namely the polarization class. By using the variational formulas from
\cite{T1} \& \cite{To89} and the facts that both $\log($det($\triangle_{\tau}%
$)) and log($\Vert\omega_{\tau}\Vert^{2})$ are potentials of the Weil
Petersson metric, we can prove that locally we have the following formula
det($\triangle_{\tau}$)$=\Vert\omega_{\tau}\Vert^{2}|f|^{2}$, where $f$ is a
holomorphic function and $\Vert\omega_{\tau}\Vert^{2}$ is a family of
\ $L^{2}$ norms of holomorphic n forms on the CY manifolds. We will prove that
$\Vert\omega_{\tau}\Vert^{2}$ has a logarithmic growth as $\tau$ approaches
$\mathcal{D}_{\infty}.$ We proved in \cite{T2} that the Ray Singer Analytic
Torsion is bounded by a constant. Combining those two facts we deduce that the
canonical C$^{\infty}$ det($\overline{\partial}$) vanishes on $\mathcal{D}%
_{\infty}.$

\textbf{Step 2. }Next step is to show that some power of the determinant line
bundle $\mathcal{L}^{N}$ is a trivial line bundle over $\mathcal{M}(M).$ This
fact follows from the following five ingredients.

\begin{enumerate}
\item  We know that as a C$^{\infty}$ line bundle $\mathcal{L}$ is trivial
over $\mathcal{M}(M).$

\item  The second fact is that $\mathcal{L}$ as a holomorphic line bundle is
isomorphic to $R^{0}\pi_{\ast}(\Omega_{\mathcal{X}/\mathcal{M}(M)}^{2n+1,0}).$
So $\mathcal{L}$ is a subbundle of the flat vector bundle $R^{n}\pi_{\ast
}\mathbb{C}.$

\item  The pull back of $R^{n}\pi_{\ast}\mathbb{C}$ is a trivial vector bundle
over the Teichm\"{u}ller space $\mathcal{T}(M)$ we can conclude that the pull
back of $\mathcal{L}^{\ast}$ on $\mathcal{T}(M)$ is also a trivial holomorphic bundle.

\item  Since $\mathcal{M}(M)\approxeq\mathcal{T}(M)/\Gamma,$ where $\Gamma$ is
some arithmetic group of rank $\geq2,$ we conclude that the holomorphic line
bundle $\mathcal{L}$ is a flat line bundle defined by character of $\Gamma.$

\item  A theorem of Kazhdan, which states that $\Gamma/[\Gamma,\Gamma]$ is
finite implies that $\mathcal{L}^{N}$ is a trivial holomorphic line bundle
over $\mathcal{M}(M).$
\end{enumerate}

From these five ingredients we deduce the existence of a non vanishing
holomorphic section $\eta$ of the determinant line bundle $\overline
{\mathcal{L}}$ over $\overline{\mathcal{M}(M)}$ \ whose zero set is supported
by the discriminant locus $\frak{D}_{\infty}.$

\subsection{The Organization of the Paper}

This article is organized as follows.

In \textbf{Section 2} we will construct the Teichm\"{u}ller space based on the
local deformation theory of CY manifolds developed in \cite{To89}.

In \textbf{Section 3 }we will prove the existence of a subgroup $\Gamma$ of
finite index in the mapping class group of a CY manifold such that this
subgroup acts freely on the Teichm\"{u}ller space. We also will show that the
quotient of the Teichm\"{u}ller space by $\Gamma$ is a non-singular one. So we
will construct a finite covering $\mathcal{M}(M)$ of the moduli space of CY
manifolds $\mathcal{M}^{\prime}(M)$, which is a non-singular quasi-projective variety.

In \textbf{Section 4} we recall the theory of determinant line bundles of
Mumford, Knudsen, Bismut, Donaldson, Gillet and Soul\'{e}, following the
exposition of D. Freed. See \cite{D} and \cite{F}. In this section we will
construct a non vanishing section $\det(D)$ of the determinant line bundle
$\mathcal{L}$ over the moduli space $\mathcal{M}$(M) of a CY manifold M of any dimension.

In \textbf{Section 5 }we prove that the determinant line bundle is a trivial
C$^{\infty}$ bundle over the moduli space $\mathcal{M}(M).$ We also construct
a unique up to a constant $\xi$ such that $|\xi|=1$ C$^{\infty}$ section
$\det(\overline{\partial})$ of the determinant line bundle $\mathcal{L}$ which
has a Quillen norm equal to Ray Singer Analytic torsion $I(M).$

In \textbf{Section 6}we review the theory of metrics with logarithmic
singularities following Mumford's article \cite{Mu}. We will show that the
determinant line bundle is isomorphic to the line bundle of holomorphic
n-forms on the moduli space $\mathcal{M}(M)$ and the natural $L^{2}$ metric on
that bundle has logarithmic singularities.

In \textbf{Section 7} we will use the results of the previous sections to
deduce that we can prolong the determinant line bundle to a line bundle on any
compactification $\overline{\mathcal{M}(M)\text{ }}$ such that $\overline
{\mathcal{M}(M)\text{ }}\backslash\mathcal{M}(M)\frak{=D}_{\infty}$ is a
divisor with normal crossings$\frak{.}$ By using the fact proved in \cite{T1}
that the analytic Ray-Singer torsion is equal to the determinant of the
Laplacian of CY metric, acting on functions, we deduce that we can prolong the
canonical C$^{\infty}$ section $\det(\overline{\partial})$ to any
compactification $\overline{\mathcal{M}(M)\text{ }}$ such that $\overline
{\mathcal{M}(M)\text{ }}\backslash\mathcal{M}(M)=\mathcal{D}_{\infty}$ is a
divisor with normal crossings and that the zero set of $\det(D)$ is supported
exactly by the discriminant locus $\frak{D}_{\infty}\frak{.}$ Based on a
result of Kazhdan and Sullivan we prove that there exists a positive integer
$N$ such that as a holomorphic line bundle the determinant line bundle to
power $N$ is a trivial one over the moduli space $\mathcal{M}(M).$ Using this
result we constructed a holomorphic section $\overline{\eta}^{N}$ of the
determinant line bundle $\left(  \overline{\mathcal{L}}\right)  ^{\otimes N}$
over $\overline{\mathcal{M}(M)}$ whose zero set is supported by $\frak{D}%
_{\infty}.$

In \textbf{Section 8} we discuss some applications of the results of this
paper and some conjectures.

\begin{acknowledgement}
I want to thank Professor Liu, Professor Li and Professor Eliashberg for
inviting me to give series of talks on the paper at Stanford University. I
want to thank Professor Yau and Professor Donaldson for their interest in the
paper and encouragement. I am deeply obliged to Professor Deligne for his
harsh and thoughtful criticism.
\end{acknowledgement}

\section{Teichm\"{u}ller Theory of CY\ Manifolds.}

\subsection{Some Definitions}

\begin{definition}
\label{Teich}We will define the Teichm\"{u}ller space $\mathcal{T}$(M) of a CY
manifold M as follows:
\end{definition}

\begin{center}
$\mathcal{T}$(M):=$\mathcal{I}$(M)/Diff$_{0}$(M),
\end{center}

\textit{where}\ $\mathcal{I}$(M)$:=\left\{  \text{the set of all integrable
complex structures on M}\right\}  $ \textit{and } Diff$_{0}$(M) \textit{is the
group of diffeomorphisms isotopic to identity. The action of the group
Diff(M}$_{0})$ \textit{is defined as follows; Let }$\phi\in$Diff$_{0}$(M)
\textit{then }$\phi$ \textit{acts on integrable complex structures on M by
pull back, i.e. if }$I\in C^{\infty}($M,$Hom($T(M),T(M)), \textit{then we
define }$\phi(I_{\Phi(\tau)})=\phi^{\ast}(I_{\Phi(\tau)})$.

We will call a pair (M; $\gamma_{1},...,\gamma_{b_{n}}$) a marked CY manifold
where M is a CY manifold and $\{\gamma_{1},...,\gamma_{b_{n}}\}$ is a basis of
$H_{n}$(M,$\mathbb{Z}$)/Tor.

\begin{remark}
\label{mark}It is easy to see that if we choose a basis of $H_{n}%
$(M,$\mathbb{Z}$)/Tor in one of the fibres of the Kuranishi family
$\mathcal{M\rightarrow K}$ then all the fibres will be marked, since as a
$C^{\infty}$ manifold $\mathcal{X}_{\mathcal{K}}\approxeq$M$\times\mathcal{K}$.
\end{remark}

\subsection{The Construction of the Teichm\"{u}ller Space}

The construction of the Teichm\"{u}ller space is based on the following lemma:

\begin{lemma}
\label{triv}Let $\pi:\mathcal{M\rightarrow K}$ be the Kuranishi family of a CY
manifold M. Let $\phi$ be a complex analytic automorphism of M such that
$\phi$ acts as an identity on $H_{n}$(M,$\mathbb{Z}$), then $\phi$ acts
trivially on $\mathcal{K}$.
\end{lemma}

\textbf{PROOF: }\ Remark \ref{mark} implies that we may suppose that the
Kuranishi family is marked. So we can define the period map
$p:\mathcal{K\rightarrow}\mathbf{P}(H^{n}($M,$\mathbb{Z)}\otimes\mathbb{C)}$
as follows:

\begin{center}
$p$((M; $\gamma_{1},..,\gamma_{b_{n}}$))$=(...,\int_{\gamma_{i}}\omega_{\tau
},...)\in\mathbf{P}(H^{n}($M,$\mathbb{Z)}\otimes\mathbb{C)}$,
\end{center}

where $\omega_{\tau}$ is a family of holomorphic n-forms over $\mathcal{K}$.
Local Torelli theorem states that $p$ is a local isomorphism. See \cite{Gr}.
So we can assume that $\mathcal{K\subset}\mathbf{P}(H^{n}($M,$\mathbb{Z)}%
\otimes\mathbb{C)}$. From here Lemma \ref{triv} follows directly.
$\blacksquare.$

\begin{theorem}
\label{teich}The Teichm\"{u}ller space $\mathcal{T}$(M) of CY manifold of
dimension n$\geq3$ exists and each connected component \ of $\mathcal{T}$(M)
is a complex analytic manifold of complex dimension $h^{n-1,1}=\dim
_{\mathbb{C}}H^{1}($M,$\Omega^{n-1})$.
\end{theorem}

\subsubsection{PROOF\ OF\ THEOREM\ \ref{teich}\textbf{ }}

We will define $\mathcal{T}$(M) as follows: Let $\frak{T}$ (M)\ be the set of
\ all marked Kuranishi families $\mathcal{X}_{\mathcal{K}}\rightarrow
\mathcal{K}$, then $\mathcal{T}$(M):=$\frak{T}$ (M)/$\symbol{126}$, where
$\symbol{126}$ is the following equivalence relation. (Notice that the points
of $\frak{T}$ (M) are pairs ($\tau,\pi^{-1}(\tau))$, where $\tau$ is a point
in some $\mathcal{K}$ and $\pi^{-1}(\tau)$ is a marked CY manifolds.) We will
say that ($\tau_{1},\pi^{-1}(\tau_{1}))\symbol{126}(\tau_{2},\pi^{-1}(\tau
_{2}))$ if and only if \ $\pi^{-1}(\tau_{1})$ and $\pi^{-1}(\tau_{2})$ are
isomorphic as marked CY manifolds. Lemma \ref{triv} together with the result
proved in \cite{To89} that the Kuranishi space is a non-singular complex
analytic space of dimension $h^{2,1}$ implies that each component of
$\mathcal{T}$(M) is a complex analytic manifold of complex dimension
$h^{n-1,1}=\dim_{\mathbb{C}}H^{1}($M,$\Omega^{n-1})$. From the
Mutsusaka-Mumford Theorem, we know that the moduli topology on $\mathcal{K}$
is a Haussdorff one. \cite{MM}. Indeed the completeness of the Kuranishi
family $\mathcal{X}_{\mathcal{K}}\rightarrow\mathcal{K}$ and the Theorem from
\cite{MM} implies that if we have two families $\mathcal{Y}_{\mathcal{K}%
}\rightarrow\mathcal{K}$ and $\mathcal{X}_{\mathcal{K}}\rightarrow\mathcal{K}$
and sequence of points$\left\{  \tau_{i}\right\}  $ such that

\begin{center}
$\underset{i\rightarrow\infty}{\lim}\tau_{i}=\tau_{0}\in\mathcal{K}$
\end{center}

such that the fibres $Y_{\tau_{i}}$ and $X_{\tau_{i}}$ are isomorphic then
$X_{\tau_{0}}$ is isomorphic to $Y_{\tau_{0}}.$ This fact is based on the
observation that the isomorphism between $Y_{\tau_{i}}$ and $X_{\tau_{i}}$
should preserve a fixed polarization. We can apply the Bishop Theorem to
conclude Theorem \ref{teich}. See \cite{Bi}. Theorem \ref{teich} is proved.
$\blacksquare.$

In Theorem \ref{teich} we proved much more. Indeed it is easy to prove that if
$\phi$ is a complex analytic automorphism of CY manifold M isotopic to
identity and $\phi$ is of a finite order then it must be the identity. So the
construction of $\mathcal{T}$(M) implies directly that we constructed an
universal family $\mathcal{M\rightarrow T}$(M) of marked CY manifolds.

From now on we will denote by $\mathcal{T}$(M) the irreducible component of
the Teichm\"{u}ller space that contains our fix CY manifold M.

\section{Construction of the Moduli Space of Polarized CY Manifolds}

The group $\Gamma_{1}:=$Diff$^{+}$(M)/Diff$_{0}$(M), where Diff$^{+}$(M) is
the group of diffeomorphisms preserving the orientation of M and Diff$_{0}$(M)
is the group of diffeomorphisms of M isotopic to identity will be called the
mapping class group. A pair (M;$L$) will be called a polarized CY manifold if
$L\in H^{2}($M,$\mathbb{Z)}$ is a fixed class and there exists a K\"{a}hler
metric G such that [ImG]$=L$. We will define $\Gamma_{2}:=\{\phi\in\Gamma
_{1}|\phi(L)=L\}.$ The reason for using the group $\Gamma_{2}$ is that the
moduli space $\Gamma_{2}$%
$\backslash$%
$\mathcal{T}$(M) will be Haussdorff.

\begin{theorem}
\label{Vie}There exists a subgroup of finite index $\Gamma$ of $\ \Gamma_{2}$
such that $\Gamma$ acts freely on $\mathcal{T}$(M) and $\Gamma$%
$\backslash$%
$\mathcal{T}$(M)$=\frak{M}$(M) is a non-singular quasi-projective variety.
\end{theorem}

\subsubsection{PROOF\ OF\ THEOREM\ \ref{Vie}}

The results of Sullivan from \cite{Sul} imply that $\Gamma_{2}$ is an
arithmetic group. This implies that in case of odd dimensional CY manifolds
there is a homomorphism induced by the action of the diffeomorphism group on
the middle homology with coefficients in $\mathbb{Z}:\phi:\Gamma
_{2}\rightarrow$Sp($2b_{n},\mathbb{Z)}$ such that the image of $\Gamma_{2}$
has a finite index in the group Sp($2b_{n},\mathbb{Z)}$\ and $\ker(\phi)$ is a
finite group.

In the case of even dimensional CY there is a homomorphism $\phi:\Gamma
_{2}\rightarrow$SO($2p,q;\mathbb{Z)}$ where SO($2p,q;\mathbb{Z})$ is the group
of the automorphisms of the lattice $H_{n}$(M,$\mathbb{Z})$/Tor, where the
$\phi(\Gamma_{2})$ has a finite index in the group SO($2p,q;\mathbb{Z)}$\ and
$\ker(\phi)$ is a finite group. A theorem of Borel implies that we can always
find a subgroup of finite index $\Gamma$ in $\Gamma_{2}$ such that $\Gamma$
acts freely on Sp($2b_{n},\mathbb{R)}$/U($b_{n}\mathbb{)}$ or on SO$_{0}%
$($2p,q;\mathbb{R})$/SO($2p$)$\times$SO($q$). We will prove that $\Gamma$ acts
without fixed point on $\mathcal{T}$(M).

Suppose that there exists an element $g\in\Gamma$, such that $g(\tau)=\tau$
for some $\tau\in\mathcal{K\subset T}$(M). From local Torelli theorem we
deduce that we may assume that the Kuranishi space $\mathcal{K}$ is embedded
in the $\mathcal{G}$, where $\mathcal{G}$ is the classifying space of Hodge
structures of weight n on $H^{n}($M,$\mathbb{Z)\otimes C}$. Griffiths proved
in \cite{Gr} that $\mathcal{G\thickapprox}G/K$ \ where G in the odd
dimensional case is Sp($2b_{n},\mathbb{R)}$ and in the even dimensional is
SO$_{0}$($2p,q;\mathbb{R})$ and $K$ is a compact subgroup of $G.$

Let $K_{0}$ is the maximal compact subgroup of $G.$ So we have a natural
$C^{\infty}$ fibration $K_{0}/K\subset G/K\rightarrow G/K_{0}$. Griffith's
transversality theorem implies that $\mathcal{K}$ is transversal to the fibres
$K_{0}/K$ of the fibration $G/K\rightarrow G/K_{0}$.

The first part of our theorem follows from the fact that $\mathcal{K}$ is
transversal to the fibres $K_{0}/K$ of the fibration $G/K\rightarrow G/K_{0}$
and the following observation; if $g\in\Gamma$ fixes a point $\tau\in G/K_{0}
$, then $g\in K_{0}\cap\Gamma$.\footnote{We suppose that $K$ \ or $K_{0}$ acts
on the right on $G$ and $\Gamma$ acts on the left on $G.$} On the other hand
side it is easy to see that local Torelli theorem implies that the action of
$\Gamma$ on $\mathcal{K}$ is induced from the action $\Gamma$ on $G/K$ by left
multiplications. From here we deduce that the action of $\Gamma$ preserve the
fibration $K_{0}/K\subset G/K\rightarrow G/K_{0}.$ From here and the fact that
$\Gamma$ acts without fix point on $G/K_{0}$ the first part of our theorem
follows directly. The second part of the theorem, namely that the space
$\Gamma$%
$\backslash$%
$\mathcal{T}$(M) is a quasi projective follows directly from the fact that
$\Gamma$%
$\backslash$%
$\mathcal{T}$(M)$\rightarrow\Gamma_{2}$%
$\backslash$%
$\mathcal{T}$(M) is a finite map and that $\Gamma_{2}$%
$\backslash$%
$\mathcal{T}$(M) is a quasi projective variety according to \cite{W}. Our
theorem is proved. $\blacksquare.$

\begin{remark}
\label{Vie1}From Theorem \ref{Vie} it follows that we constructed a family of
non-singular CY manifolds $\mathcal{X\rightarrow M}$(M) over a
quasi-projective non-singular variety $\mathcal{M}$(M). Moreover it is easy to
see that $\mathcal{X\subset}\mathbb{CP}^{N}\times\mathcal{M}$(M). \textit{So}
$\mathcal{X}$ \ \textit{is also quasi-projective. From now on we will work
only with this family.}
\end{remark}

\section{The Theory of Determinant Line Bundles}

\subsection{Geometric Data}

In order to construct the determinant line bundle we need the following data:

\begin{enumerate}
\item  A smooth fibration of manifolds $\pi:\mathcal{X\rightarrow M(}$M). In
our case it will be the smooth fibration of the family of CY manifolds over
the moduli space as defined in Theorem \ref{Vie}. Let n=$\dim_{\mathbb{C}}M.$

\item  A metric along the fibres, that is a metric g($\tau)$ on the relative
tangent bundle $\mathcal{T(X}$/$\mathcal{M}$(M)). In this paper the metric
will be the families of CY metrics g($\tau$) such that the class of cohomology
$[\operatorname{Im}(g(\tau))]=L$ \ is fixed.
\end{enumerate}

From these data we will construct the determinant line bundle $\mathcal{L}$
over the moduli space of $\ $CY manifolds $\mathcal{M}$(M). We will consider
the relative $\overline{\partial}_{\mathcal{X}\text{/}\mathcal{M}\text{(M)}}$ complex:

\begin{center}
$0\rightarrow$ker$\overline{\partial}_{\mathcal{X}\text{/}\mathcal{M}%
\text{(M)}}\rightarrow$C$^{\infty}(\mathcal{X}$)$\overset{\overline{\partial
}_{0,\mathcal{X}\text{/}\mathcal{M}\text{(M)}}}{\rightarrow}\Omega
_{\mathcal{X}\text{/}\mathcal{M}\text{(M)}}^{0,1}\overset{\overline{\partial
}_{1,\mathcal{X}\text{/}\mathcal{M}\text{(M)}}}{\rightarrow}..\Omega
_{\mathcal{X}\text{/}\mathcal{M}\text{(M)}}^{0.n-1}\overset{\overline
{\partial}_{n-1,\mathcal{X}\text{/}\mathcal{M}\text{(M)}}}{\rightarrow}%
.\Omega_{\mathcal{X}\text{/}\mathcal{M}\text{(M)}}^{0.n}\rightarrow0.$
\end{center}

We will define $D$ to be for each $\tau\in\mathcal{M}$(M) and k,

\begin{center}
$D_{k}:=\overline{\partial}_{k,\mathcal{X}\text{/}\mathcal{M}\text{(M)}%
}+\left(  \overline{\partial}_{k,\mathcal{X}\text{/}\mathcal{M}\text{(M)}%
}\right)  ^{\ast}\&$ $D_{k,\tau}:=D_{k}\left|  _{M_{\tau}}\right.
=\overline{\partial}_{k,\tau}+\left(  \overline{\partial}_{k,\tau}\right)
^{\ast}.$
\end{center}

\begin{definition}
\label{DR}We will call the above complex the\ relative Dolbault complex.
\end{definition}

Let us define $\left(  \mathcal{H}^{k}\right)  _{\tau}:=L^{2}($M$_{\tau
},\Omega_{\tau}^{0,k}).$ Furthermore, as $\tau$ varies over $\mathcal{M}$(M),
the spaces $\left(  \mathcal{H}_{\tau}^{k}\right)  $ fit together to form
continuous Hilbert bundles $\mathcal{H}^{k}$ over $\mathcal{M}$%
(M).\footnote{These bundles are not smooth since the composition $L^{2}\times
C^{\infty}\rightarrow L^{2}$ is not differentiable map.}. Thus we can view
$\overline{\partial}_{k,\mathcal{X}\text{/}\mathcal{M}\text{(M)}}$ as a bundle
maps: $\overline{\partial}_{k,\mathcal{X}\text{/}\mathcal{M}\text{(M)}%
}:\mathcal{H}^{k}\rightarrow\mathcal{H}^{k+1}.$ The Hilbert bundles
$\mathcal{H}^{k}$ carry $L^{2}$ metrics by definition.

\subsection{Construction of the Determinant Line Bundle $\mathcal{L}$}

\subsubsection{Some Basic Definitions}

Now we are ready to construct the Determinant line bundle $\mathcal{L}$ of the
operator $\overline{\partial}_{\mathcal{X}\text{/}\mathcal{M}\text{(M)}}.$ We
will recall some basic consequences of the ellipticity if $D_{\tau}$. Each
fibre $\mathcal{H}_{\tau}^{k}$ of the Hilbert bundles $\mathcal{H}^{k}$
decomposes into direct sum of eigen spaces of non-negative Laplcaians
$D_{k}D_{k}^{\ast}$ and $D_{k}^{\ast}D_{k}.$ The spectrums of these operators
are discrete, and the nonzero eigen values $\{\lambda_{k,i}\}$ of
\ $D_{k}^{\ast}D_{k}$ and $D_{k}D_{k}^{\ast}$ agree and D$_{k}$ define a
canonical isomorphisms between the corresponding eigen spaces.

\begin{definition}
\label{1}\textbf{1.} Let $\mathcal{U}_{a}:\{\tau\in\mathcal{M}(M)|a\notin
Spec(D_{k}D_{k}^{\ast})$ for $0\leq k\leq n$ and any $a>0\}.$ \textbf{(}%
$U_{a}$ are open sets in $\mathcal{M}(M)$ and they form an open covering of
$\mathcal{M}$(M) since the spectrum of $D_{\tau}^{\ast}D_{\tau}$ is discrete.)
\textbf{2. }Let the fibres of $\mathcal{H}_{a}^{k}$ be the vector subspaces in
$\mathcal{H}_{\tau,a}^{k}$ spanned by eigen vectors with eigen values less
than $a$ over $\mathcal{U}_{a}.$ Then we can define the complex of :
\end{definition}

\begin{center}
$0\rightarrow\Gamma(\mathcal{U}_{a},\mathcal{O}_{\mathcal{U}_{a}}%
)\rightarrow\mathcal{H}_{a}^{0}\overset{\overline{\partial}_{0,\mathcal{X}%
\text{/}\mathcal{M}\text{(M)}}}{\rightarrow}...\rightarrow\mathcal{H}%
_{a}^{n-1}\overset{\overline{\partial}_{n-1,\mathcal{X}\text{/}\mathcal{M}%
\text{(M)}}}{\rightarrow}\mathcal{H}_{a}^{n}\rightarrow\ker(D_{n-1}\circ
D_{n}^{\ast})\rightarrow0.$
\end{center}

If $b>a$ we set $\mathcal{H}_{a,b}^{k}:=\mathcal{H}_{b}^{k}$/$\mathcal{H}%
_{a}^{k}.$ The spaces $\mathcal{H}_{a}^{k}$ form a smooth finite dimensional
$C^{\infty}$ bundles over an open set $\mathcal{U}^{a}\subset\mathcal{M}$(M).
For the proof of this fact see \cite{ABKS}.

\subsubsection{Construction of the Generating Sections $\det(D_{a})$ over
$\mathcal{U}_{a}$}

\begin{definition}
\label{3} Let $\omega_{1}^{k},..,\omega_{m_{k}}^{k},\psi_{1}^{k}%
,..,\psi_{N_{k}}^{k},\phi_{1}^{k},..,\phi_{M_{k}}^{k}$ be an orthonormal basis
in the trivial vector bundle $\mathcal{H}_{a}^{k}$ over $\mathcal{U}_{a}$,
where $D_{k}\omega_{i}^{k}=0$, $\overline{\partial}_{k}^{\ast}(\overline
{\partial}_{k}\psi_{j}^{k})=\lambda_{j}^{k}\psi_{j}^{k},$ $\overline{\partial
}_{k}(\overline{\partial}_{k}^{\ast}\phi_{j}^{k})=\lambda_{j}^{k}\phi_{j}%
^{k},$ $\phi_{j}^{k}\in\operatorname{Im}\overline{\partial}_{k-1,\mathcal{X}%
\text{/}\mathcal{M}\text{(M)}}$ and $\psi_{j}^{k}\in\operatorname{Im}%
(\overline{\partial}_{k,\mathcal{X}\text{/}\mathcal{M}\text{(M)}}^{\ast}%
)$\ for 1$\leq i\leq k$ and $0<\lambda_{j}<a$ for 1$\leq j\leq N$. Let
\end{definition}

\begin{center}
$\det(\overline{\partial}_{k,a})=\left(  \omega_{1}^{k}\wedge...\wedge
\omega_{m_{k}}^{k}\wedge(\overline{\partial}_{k-1,\mathcal{X}\text{/}%
\mathcal{M}\text{(M)}}\psi_{1}^{k-1})\wedge...\wedge(\overline{\partial
}_{k-1,\mathcal{X}\text{/}\mathcal{M}\text{(M)}}\psi_{N_{k}}^{k-1})\wedge
\phi_{1}^{k}\wedge...\wedge\phi_{M_{k}}^{k}\right)  ^{(-1)^{k}}.$
\end{center}

We will define the line bundle $\mathcal{L}$ restricted on $\mathcal{U}_{a}$
as follows:

\begin{center}
$\mathcal{L}^{a}:=\mathcal{L}|_{\mathcal{U}_{a}}=\otimes_{k=0}^{n}\left(
\wedge^{\dim\mathcal{H}_{a}^{k}}\mathcal{H}_{a}^{k}\right)  ^{(-1)^{k}}.$
\end{center}

\begin{definition}
\label{2}The generator $\det(\overline{\partial}_{a})$ of $\mathcal{L}%
^{a}:=\mathcal{L}|_{\mathcal{U}^{a}}$ is defined as follows: $\det
(\overline{\partial}_{a}):=\otimes\det(\overline{\partial}_{k,a}).$
\end{definition}

\subsubsection{Definition of the Transition Functions on $\mathcal{U}_{a}%
\cap\mathcal{U}_{b}$}

We will define how we patch together $\mathcal{L}^{a}$ and $\mathcal{L}^{b}$
over $\mathcal{U}^{a}\cap\mathcal{U}^{b}\footnote{We may suppose that $b>a.$%
}.$ On that intersection we have: $\mathcal{L}^{b}=\mathcal{L}^{a}%
\otimes\mathcal{L}^{a,b},$ where$\mathcal{L}^{a,b}:=\otimes_k=0^{n}%
(\det\mathcal{H}_a,b^{k})^{(-1)^{k}}$ on $\ \mathcal{U}^{a}\cap\mathcal{U}%
^{b}.$ We can identify $\mathcal{L}^{a,b}$ over \ $\mathcal{U}^{a}%
\cap\mathcal{U}^{b}$ with the line bundle bundle spanned by the section

\begin{center}
$\det(\overline{\partial}_{a,b})=\otimes_{k=0}^{n}\det(\overline{\partial
}_{k,a,b})^{(-1)^{k}},$
\end{center}

where $\det(\overline{\partial}_{k,a,b}):=\left(  (\overline{\partial
}_{k-1,\mathcal{X}\text{/}\mathcal{M}\text{(M)}}\psi_{1}^{k-1})\wedge
...\wedge(\overline{\partial}_{k-1,\mathcal{X}\text{/}\mathcal{M}\text{(M)}%
}\psi_{N_{k}}^{k-1})\wedge\phi_{1}^{k}\wedge...\wedge\phi_{M_{k}}^{k}\right)
,$ $\phi_{j}^{k}\in\operatorname{Im}\overline{\partial}_{k-1,\mathcal{X}%
\text{/}\mathcal{M}\text{(M)}}$, $\psi_{j}^{k}\in\operatorname{Im}%
(\overline{\partial}_{0,\mathcal{X}\text{/}\mathcal{M}\text{(M)}}^{\ast})$,
$\Delta_{k}\phi_{j}^{k}=\lambda_{j}^{k}\phi_{j}^{k},$ $\Delta_{k}%
(\overline{\partial}(\psi_{i}^{k-1}))=\lambda_{i}^{k}(\overline{\partial}%
(\psi_{i}^{k-1}))$ and $a<\lambda_{i}^{k}<b.$

\begin{remark}
\label{ciso}We can view $\det(\overline{\partial}_{a,b})$ as a section of the
line bundle $\mathcal{L}^{a,b}$ over $\mathcal{U}^{a}\cap\mathcal{U}^{b}$ and
it defines a canonical smooth isomorphisms over $\mathcal{U}^{a}%
\cap\mathcal{U}^{b}:\mathcal{L}^{a}\rightarrow\mathcal{L}^{a}\otimes
\mathcal{L}^{a,b}=\mathcal{L}^{b}$ $(s\rightarrow s\otimes\det(\overline
{\partial}_{a,b}).$
\end{remark}

We define the determinant line bundle $\mathcal{L}$ by patching the trivial
line bundles $\mathcal{L}^{a}$ over $\mathcal{U}^{a}$ by using the canonical
isomorphism defined in Remark \ref{ciso}.

\subsection{The Description of the Quillen Metric on $\mathcal{L}$}

We now proceed to describe the Quillen metric on $\mathcal{L}$. Fix $a>0$.
Then the subbundles $\mathcal{H}_{a}^{k}$ of the Hilbert bundles
$\mathcal{H}^{k}$ on $\mathcal{U}_{a}$ inherit metrics from $\mathcal{H}^{k}.$
According to standart facts from linear algebra, metrics are induced on
determinants, duals, and tensor products. So the $\mathcal{L}^{a}$ inherits a
natural metric. We will denote by g$^{a}$ the $L^{2}$ norm of the section
$\det(\overline{\partial}_{a})$. Clearly

\begin{center}
g$^{a}=\prod_{k=0}^{n}\left(  \lambda_{1}^{k}...\lambda_{n_{k}}^{k}\right)
^{(-1)^{k}},$
\end{center}

where the product is of all non zero eigen values of the operators
$\overline{\partial}_{k}^{\ast}\overline{\partial}_{k-1}$which are less than $a.$

If $b>a$, then under the isomorphism defined in Remark \ref{ciso}, we have two
metrics on $\mathcal{L}^{b}$ and their ratio is a real number equal to the
$L^{2}$ norm of the section$\left\|  \det(\overline{\partial}_{a,b})\right\|
^{2}$. The definition of the section $\det(\overline{\partial}_{a,b})$ implies
that we have the following formula:

\begin{center}
$\left\|  \det(\overline{\partial}_{a,b})\right\|  =\prod_{k=0}^{n}\prod
_{i=1}\left\|  \phi_{i}^{k}\right\|  \prod_{j=1}\left\|  \overline{\partial
}\psi_{j}^{k}\right\|  ^{(-1)^{k}}=\prod_{i=1}\left\|  \phi_{i}^{k}\right\|
\prod_{j=1}\left\langle \overline{\partial}_{k}^{\ast}\overline{\partial
}_{k-1}\psi_{j}^{k},\psi_{j}^{k}\right\rangle ^{(-1)^{k}}=\prod\left(
\lambda_{i}^{k}\right)  ^{(-1)^{k}}.$
\end{center}

where $\lambda_{i}^{k}$ are all the non-zero eigen values of the operators
$\overline{\partial}_{k}^{\ast}\overline{\partial}_{k-1}$ such that
$a<\lambda_{i}^{k}<b.$ In other word, on $\mathcal{U}^{a}\cap\mathcal{U}^{b}$
g$^{b}=$g$^{a}\prod\left(  \lambda_{i}^{k}\right)  ^{(-1)^{k}}.$ To correct
this discrepancy we define $\overline{g}^{a}=g^{a}\det(D^{\ast}D\left|
_{\lambda>a}\right.  ),$ where

\begin{center}
$\det(\overline{\partial}_{k}^{\ast}\overline{\partial}_{k-1}\left|
_{\lambda>a}\right.  )=-\exp(-\left(  \zeta_{k}^{a}\right)  ^{^{\prime}}(0))$
and $\zeta_{k}^{a}(s)=\sum_{\lambda_{i}>a}^{\infty}\left(  \lambda_{i}%
^{k}\right)  ^{s}.$
\end{center}

The crucial property of this regularization is that it behaves properly with
respect to the finite number of eigen values, i.e.

\begin{center}
$\det(\overline{\partial}_{k}^{\ast}\overline{\partial}_{k-1}\left|
_{\lambda>b}\right.  )=\det(\overline{\partial}_{k}^{\ast}\overline{\partial
}_{k-1}\left|  _{\lambda>a}\right.  )\prod_{a<i<b}^{N}\lambda_{i}^{k}$
\end{center}

on the intersection $\mathcal{U}^{a}\cap\mathcal{U}^{b}.$ From the last remark
we deduce that $\overline{g}^{a}$ and $\overline{g}^{b}$ agree on
$\mathcal{U}^{a}\cap\mathcal{U}^{b}.$ Thus $\overline{g}^{a}$ patch together
to a Hermitian metric g$^{\mathcal{L}}$ on $\mathcal{L}$. The metric
g$^{\mathcal{L}}$ will be called the Quillen metric on $\mathcal{L}$.

\begin{definition}
\label{rs}We will define the holomorphic analytic torsion I(M), for odd
dimensional CY manifold M as follows: I(M)$:=\prod_{q=1}^{n}(\det
(\triangle_{q}^{^{\prime}})^{(-1)^{q}}.$ See \cite{RS}.
\end{definition}

\begin{remark}
\label{rs1}It is easy to see that if dim$_{\mathbb{C}}$M=2n, then $\log
$I(M)=0. We know from the results in \cite{T1} that for odd dimensional CY
manifolds $I(M)>0.$ So from now on we will consider only odd dimensional CY manifolds.
\end{remark}

We will need the following result from \cite{BGS1} on p. 55:

\begin{theorem}
\label{BGS2}The Quillen norm of the $C^{\infty}$ section $\det(\overline
{\partial}_{a})$ on $\mathcal{U}_{a}$ of $\mathcal{L}$ is equal to I(M).
\end{theorem}

\subsubsection{PROOF\ OF\ THEOREM\ \ref{BGS2}}

It follows from the Definition \ref{3}\ of the section $\det(\overline
{\partial})\left|  _{\mathcal{U}^{a}}\right.  $ of $\mathcal{L}$ and the
definition of the Quillen metric that at each point $\tau\in\mathcal{M}$(M)
the following formula is true:$\left\|  \det(\overline{\partial})_{\tau
}\left|  _{\mathcal{U}^{a}}\right.  \right\|  _{Q}^{2}=I(M_{\tau
})|_{\mathcal{U}^{a}},$ where $\left\|  \det(\overline{\partial})_{\tau
}\left|  _{\mathcal{U}^{a}}\right.  \right\|  _{Q}^{2}$ means the Quillen norm
of the section $\det(\overline{\partial})_{\tau}\left|  _{\mathcal{U}^{a}%
}\right.  .$ Theorem \ref{BGS2} is proved. $\blacksquare.$

\section{Construction of a C$^{\infty}$ Non Vanishing Section of the
Determinant Line Bundle $\mathcal{L}$ for Odd Dimensional CY Manifolds}

\subsection{Some Ppreliminary Results}

Let us denote by $\pi_{\ast}\left(  \omega_{\mathcal{X}\text{/}\mathcal{M}%
\text{(M)}}\right)  :=\pi_{\ast}\left(  \Omega_{\mathcal{X}\text{/}%
\mathcal{M}\text{(M)}}^{n,0}\right)  $ the relative dualizing sheaf. The local
sections of $\pi_{\ast}\left(  \omega_{\mathcal{X}\text{/}\mathcal{M}%
\text{(M)}}\right)  $ are families of holomorphic n-forms $\omega_{\tau}$ on
M$_{\tau}$. Then on $\pi_{\ast}\left(  \omega_{\mathcal{X}\text{/}%
\mathcal{M}\text{(M)}}\right)  $ we have a natural $L^{2}$ metric defined as follows:

\begin{center}
$\left\|  \omega_{\tau}\right\|  ^{2}:=(-1)^{\frac{n(n-1)}{2}}\left(
\sqrt{-1}\right)  ^{n}\int_{\text{M}}\omega_{\tau}\wedge\overline{\omega
_{\tau}}. $
\end{center}

\begin{theorem}
\label{Tod1}If the dimension of the CY manifold is even then $\mathcal{L}$
$\ $is isomorphic to the line bundle $\pi_{\ast}\left(  \omega_{\mathcal{X}%
\text{/}\mathcal{M}\text{(M)}}\right)  .$ If the dimension of the CY manifold
is odd then the dual of $\mathcal{L}$ is isomorphic to the line bundle
$\pi_{\ast}\left(  \omega_{\mathcal{X}\text{/}\mathcal{M}\text{(M)}}\right)  $
over $\mathcal{M}(M).$
\end{theorem}

\subsubsection{PROOF\ OF\ THEOREM\ \ref{Tod1}}

From the definition of CY says that:

\begin{center}
$\dim_{\mathbb{C}}H^{j}($M,$\mathcal{O}_{\text{M}}$)=$\left\{
\begin{array}
[c]{ll}%
1 & j=0\text{ or }j=n\\
0 & for\text{ }j\neq0\text{ or }n
\end{array}
\right.  $
\end{center}

This and the definition of CY manifold imply that

\begin{center}
\bigskip$R^{q}\pi_{\ast}\mathcal{O}_{\text{M}}=\left\{
\begin{array}
[c]{ll}%
\left(  \pi_{\ast}\omega_{\mathcal{X}\text{/}\mathcal{M}\text{(M)}}\right)
^{\ast} & \text{ }j=n\\
\mathcal{O}_{\mathcal{M}\text{(M)}} & for\text{ }j\neq\text{ }n
\end{array}
\right.  $
\end{center}

From the definition of $\mathcal{L}$ it follows that $\mathcal{L\backsimeq
}\prod_{q=0}^{n}\left(  -1\right)  ^{q}\det\left(  R^{q}\pi_{\ast}%
\mathcal{O}_{\text{M}}\right)  .$ Combining all these equalities we deduce
directly Theorem \ref{Tod1}. $\blacksquare.$

\begin{corollary}
\label{zero}Let M be a CY manifold of odd dimension $n=2m+1.$ Then the index
of the operator $\overline{\partial}$ on the complex defined in Definition
\ref{DR} is zero.
\end{corollary}

\subsection{Holomorphic Structure on the Determinant line Bundle $\mathcal{L}$}

In \cite{BGS1} a canonical smooth isomorphism is constructed between the
holomorphic determinant of Grothendieck-Knudsen-Mumford with Quillen
determinant bundle. More precisely the following theorem is proved:

\begin{theorem}
\label{BGS1}Let $\pi:\mathcal{X\rightarrow M}$(M) be a holomorphic fibration
with smooth fibres. Suppose $\mathcal{X}$ admits a closed (1,1) form $\psi$
which restricts to a K\"{a}hler metric on each fibre. Let
$\mathcal{E\rightarrow X}$ \ be a holomorphic Hermitian bundle with its
Hermitian connection. Then the determinant line bundle $\mathcal{L\rightarrow
M}$(M) of the relative $\overline{\partial}$ complex (coupled to $\mathcal{E}$
) admits a holomorphic structure. The canonical connection (constructed in
\cite{BF}) on $\mathcal{L}$ is the Hermitian connection for the Quillen
metric. Finally, if the index of \ $\overline{\partial}$ is zero, the section
$\det($ $\overline{\partial}_{E})$ of $\mathcal{L}$ is holomorphic.
\end{theorem}

From now on we will consider the family of CY manifolds $\mathcal{X\rightarrow
M}$(M) as defined in Theorem \ref{Vie1} together with the trivial line bundle
$\mathcal{E}$ \ over $\mathcal{M}$(M). It is easy to see that the family
$\mathcal{X\rightarrow M}$(M) fulfill the conditions of the Theorem \ref{BGS1}.

\subsection{Construction of a C$^{\infty}$ Section of the Determinant Line
Bundle $\mathcal{L}$ with Quillen Norm Ray-Singer Equal to Ray-Singer Analytic Torsion}

\begin{definition}
\label{4}Let $\mathcal{H}_{+}=\underset{k}{\oplus}L^{2}(M,\Omega_{M}^{0,2k})$
and $\mathcal{H}_{-}=\underset{k}{\oplus}L^{2}(M,\Omega_{M}^{0,2k+1})$ and
$D=\overline{\partial}_{\mathcal{X}/\mathcal{M}(M)}+\overline{\partial
}_{\mathcal{X}/\mathcal{M}(M)}^{\ast}.$
\end{definition}

\begin{theorem}
\label{CYsec}Let M be a CY manifold of odd dimension, then as C$^{\infty}$
bundle the determinant line bundle $\mathcal{L}$ is trivial and there exists a
global $C^{\infty}$ section $\det(\overline{\partial})$ of
$\mathcal{L\rightarrow M}$(M), which has no zeroes on $\mathcal{M}$(M) and
whose Quillen norm is the Ray Singer Analytic Torsion.
\end{theorem}

\subsubsection{PROOF\ OF\ THEOREM\ \ref{CYsec}}

\textbf{PROOF: }The proof of Theorem \ref{CYsec} is based on the following
three Theorems:

\begin{theorem}
\label{var1}The first Chern class of the ralative dualizing sheaf $\pi_{\ast
}\omega_{\mathcal{X}/\mathcal{M}(M)}$ is given locally by the formula
$c_{1}(\pi_{\ast}\omega_{\mathcal{X}/\mathcal{M}(M)})=dd^{c}\left(
\left\langle \omega_{\tau},\omega_{\tau}\right\rangle \right)
=-\operatorname{Im}\left(  \text{Weil-Petersson metric}\right)  ,$ where
$\omega_{\tau}$ is a holomorphic family of holomorphic n forms$..$ (See
\cite{To89}.)
\end{theorem}

\begin{theorem}
\label{var}Locally on the moduli space $\mathcal{M}(M)$ the following formula holds
\end{theorem}

\begin{center}
$dd^{c}\left(  \log(\det(\triangle_{0}))\right)  =-\operatorname{Im}\left(
\text{\textit{Weil-Petersson metric}}\right)  .$ (See \cite{T1}.)
\end{center}

\begin{theorem}
\label{LY}\textbf{i.} I(M)=$\left(  \det(\triangle_{0})\right)  ^{2}$
\textbf{ii.} around each point $\tau\subset\mathcal{U}_{\tau}\in
\mathcal{M}(M)$ we have I(M)%
$\vert$%
$_{\mathcal{U}_{\tau}}=<\omega,\omega>|\phi|^{2},$ where $\phi$ is a
holomorphic function on $\mathcal{U}$ and \textbf{iii. }The positive function
$\det(\triangle_{0})$ is bounded by a constant $C,$ i.e. $\det(\triangle
_{0})<C.$ (See \cite{T1}.)
\end{theorem}

We will proof the following Lemma:

\begin{lemma}
\label{ciso3}The first Chern class $c_{1}(\mathcal{L})$ of the C$^{\infty}$
determinant line bundle $\mathcal{L}$ is equal to zero in $H^{2}%
(\mathcal{M}(M),\mathbb{Z}).$
\end{lemma}

\textbf{Proof of Lemma }\ref{ciso3}: Theorem \ref{Tod1} implies that when the
CY manifold $M$ has an odd dimension then the determinant line bundle
$\mathcal{L}$ is isomorphic to the relative dualizing sheaf $\pi_{\ast}%
(\omega_{\mathcal{X}/\mathcal{M}(M)}).$ So we need to proof that $c_{1}%
(\pi_{\ast}(\omega_{\mathcal{X}/\mathcal{M}(M)}))=0.$ The proof of Lemma
\ref{ciso3} is based on the follwoing obeservation: Notice that the definition
of the Ray Singer analytic torsion implies that $I(M_{\tau})$ is a positive
function different from a constant on $\mathcal{M}(M).$ From Theorem \ref{LY}
we know that we have the following local expression of $I(M_{\tau}):$

\begin{center}
$I(M_{\tau})|_{\mathcal{U}}=\left\|  \omega_{\tau}\right\|  ^{2}|\phi|^{2}$
\end{center}

where $\omega_{\tau}$ is a holomorphic family of holomorphic n-forms on
$M_{\tau},$ $\phi$ is a holomorphic function on $\mathcal{U}$ and

\begin{center}
$\left\|  \omega_{\tau}\right\|  ^{2}=\left(  \sqrt[2]{-1}\right)
^{n}(-1)^{\frac{n(n-1)}{2}}\int_{M_{\tau}}\omega_{\tau}\wedge\overline
{\omega_{\tau}}.$
\end{center}

Theorems \ref{var} and \ref{var1} imply that

\begin{center}
$\frac{\sqrt[2]{-1}}{2}\partial\overline{\partial}\log(I(M_{\tau}%
))=c_{1}(\mathcal{L})=c_{1}(\pi_{\ast}(\omega_{\mathcal{X}/\mathcal{M}(M)})).$
\end{center}

Since $\partial\overline{\partial}\log(I(M_{\tau}))=d\left(  \overline
{\partial}\log(I(M_{\tau}))\right)  $ we deduce that

\begin{center}
$\frac{\sqrt[2]{-1}}{2}\partial\overline{\partial}\log(I(M_{\tau}%
))=c_{1}(\mathcal{L})=\frac{\sqrt[2]{-1}}{2}d\left(  \overline{\partial}%
\log(I(M_{\tau}))\right)  =d\alpha.$
\end{center}

So $c_{1}(\mathcal{L})=0$ in $H^{2}(\mathcal{M}(M),\mathbb{Z}).$ This proves
Lemma \ref{ciso3}. $\blacksquare.$

\begin{corollary}
\label{ciso31}The determinant line bundle $\mathcal{L}$ as a $C^{\infty}$
bundle is trivial.
\end{corollary}

So the first part of Theorem \ref{CYsec} is proved. The proof of the second
part of Theorem \ref{CYsec} is based on the following Lemma:

\begin{lemma}
\label{ciso4}There exists a non vanishing global section $\det(\overline
{\partial})$ of $\ $the dual of the determinant line bundle $\mathcal{L}$ such
that the Quillen norm of det$(\overline{\partial})=I(M).$
\end{lemma}

\textbf{Proof of Lemma }\ref{ciso4}: From Corollary \ref{ciso31} we can
conclude the existence of a global $C^{\infty}$ section $\omega$ of
$\mathcal{L\rightarrow M}$(M), which has no zeroes on $\mathcal{M}$(M) and
such that for each $\tau\in\mathcal{M}(M)$ it has $L^{2}$ norms 1, i.e. we
have $\left\|  \omega_{\tau}\right\|  ^{2}=1.$ \ Since $M_{\tau}$ is an odd
dimensional CY manifold we know from Theorem \ref{Tod1} hat $\mathcal{L}$ is
isomorphic to $\pi_{\ast}(\omega_{\mathcal{X}/\mathcal{M}(M)}).$ The non
vanishing section $\omega$ of the determinant line bundle $\mathcal{L}$ can be
interpreted as a family of (2n+1,0) forms $\omega$ which generate the kernel
of $D^{\ast}:\mathcal{H}_{-}\rightarrow\mathcal{H}_{+}.$ The kernal of
$D:\mathcal{H}_{+}\rightarrow\mathcal{H}_{-}$ is generated by the constant 1.
This follows directly from the definition of the CY manifold.

From Definition \ref{3} of the section $\det(\overline{\partial}_{a})$ on the
open set $\mathcal{U}_{a}\subset\mathcal{M}(M)$, the existence of a
C$^{\infty}$ family of antiholomoprhic forms $\omega_{\tau}$ with $L^{2}$ norm
1, which trivializes $R^{2n+1}\pi_{\ast}(\mathcal{O}_{\mathcal{X}}$ $)$ over
$\mathcal{M}(M)$ and the definition of the transition functions $\{\sigma
_{a,b}\}$\ of $\mathcal{L}$ with respect to the covering $\{\mathcal{U}%
_{a}\},$ we deduce that for $b>a$ we have on $\mathcal{U}_{a}\cap
\mathcal{U}_{b}$ $\det(\overline{\partial}_{b})=$ $\det(\overline{\partial
}_{b})\sigma_{a,b}.$ This fact and Theorem \ref{BGS2} imply that we
constructed a global $C^{\infty}$ section $\det(\overline{\partial})$ of
$\mathcal{L}$ whose Quillen norm is the Ray Singer Analytic Torsion. So the
determinant line bundle $\mathcal{L}$ is a trivial $C^{\infty}$ line bundle.
Theorem \ref{CYsec} is proved. $\blacksquare.$

\begin{corollary}
\label{LY3}The determinant line bundle as a holomorphic bundle is flat over
$\mathcal{M}(M).$
\end{corollary}

\section{The Analogue of the Dedekind Eta Function for Odd Dimensional CY Manifolds}

\subsection{Metrics on Vector Bundles with Logarithmic Growth}

In Theorem \ref{Vie} we constructed the moduli space $\mathcal{M}$(M) of CY
manifolds. From the results in \cite{W} and Theorem \ref{Vie} we know that
$\mathcal{M}$(M) is a quasi-projective non-singular variety. Using Hironaka's
resolution theorem, we may suppose that $\mathcal{M}$(M)$\subset
\overline{\mathcal{M}\left(  \text{M}\right)  },$ where$\overline
{\mathcal{M}\left(  \text{M}\right)  }\backslash\mathcal{M}\left(
\text{M}\right)  =\mathcal{D}_{\infty}$is a divisor with normal crossings. We
need to show, now how we will extend the determinant line bundle $\mathcal{L}$
to a line bundle $\overline{\mathcal{L}\text{ }}$ to $\overline{\mathcal{M}%
\left(  \text{M}\right)  }.$ For this reason we are going to recall the
following definitions and results from \cite{Mu}. We will look at the
policylinders D$^{N}\subset\overline{\mathcal{M}\left(  \text{M}\right)  },$
where D is the unit disk and $N=\dim\overline{\mathcal{M}\left(
\text{M}\right)  }$ and

\begin{center}
D$^{N}\cap\mathcal{D}_{\infty}=\{$union of hyperplanes; $\tau_{1}%
=0,...,\tau_{k}=0\}.$
\end{center}

Hence D$^{N}\cap\mathcal{M}$(M)=(D$^{\ast})^{k}\times$D$^{N\leq-k}.$ In
D$^{\ast}$ we have the Poincare metric

\begin{center}
ds$^{2}=\frac{\left|  dz\right|  ^{2}}{\left|  z\right|  ^{2}\left(
\log\left|  z\right|  \right)  ^{2}}$
\end{center}

and in D we have the simple metric $\left|  dz\right|  ^{2},$ giving us a
product metric on (D$^{\ast}$)$^{k}\times$D$^{r-k}$ which we call
$\omega^{(P)}.$

A complex-valued C$^{\infty}$ p-form $\eta$ on $\mathcal{M}$(M) is said to
have Poincare growth on $\overline{\mathcal{M}\left(  \text{M}\right)
}\backslash\mathcal{M}\left(  \text{M}\right)  $ if there is a set of
\textit{if policylinders \ }$\mathcal{U}_{\alpha}\subset\overline
{\mathcal{M}\left(  \text{M}\right)  }$ covering $\overline{\mathcal{M}\left(
\text{M}\right)  }\backslash\mathcal{M}\left(  \text{M}\right)  $ such that in
each \textit{\ }$\mathcal{U}_{\alpha}$ \textit{an estimate of the following
type holds:}

\begin{center}
$\left|  \eta(\tau_{1},...,\tau_{N}\right|  \leq C_{\alpha}\omega
_{\mathcal{U}_{\alpha}}^{(P)}(\tau_{1},\overline{\tau_{1}})..._{\alpha}%
\omega_{\mathcal{U}_{\alpha}}^{(P)}(\tau_{N},\overline{\tau_{N}}).$
\end{center}

This property is independent of the covering $\mathcal{U}_{\alpha}$ of
$\overline{\mathcal{M}\left(  \text{M}\right)  }\backslash\mathcal{M}\left(
\text{M}\right)  $ but depends on the compactification $\overline
{\mathcal{M}\left(  \text{M}\right)  }.$ If $\eta_{1}$ and $\eta_{2}$ both
have Poincare growth on $\overline{\mathcal{M}\left(  \text{M}\right)
}\backslash\mathcal{M}\left(  \text{M}\right)  ,$ then so does $\eta_{1}%
\wedge\eta_{2}$. The basic property of the Poincare growth is the following:

\begin{theorem}
\label{Mum}A p-form $\eta$ with a Poincare growth on $\overline{\mathcal{M}%
\left(  \text{M}\right)  }\backslash\mathcal{M}\left(  \text{M}\right)  ,$
\textit{has the property that for every C}$^{\infty}$ (r-p) \textit{form}
$\psi$ \textit{on \ }$\overline{\mathcal{M}\left(  \text{M}\right)  }$
\textit{we have:}
\end{theorem}

\begin{center}
$\int_{\overline{\mathcal{M}\left(  \text{M}\right)  }\backslash
\mathcal{M}\left(  \text{M}\right)  }\left|  \eta\wedge\psi\right|  <\infty.$
\end{center}

\textit{Hence, }$\eta$ \textit{defines a current [}$\eta$] \textit{on
}$\overline{\mathcal{M}\left(  \text{M}\right)  }.$

\subsubsection{PROOF\ OF\ THEOREM\ \ref{Mum}}

For the proof see \cite{Mu}. $\blacksquare.$

A complex valued C$^{\infty}$ p-form\textit{\ }$\eta$ on\textit{\ }%
$\overline{\mathcal{M}\left(  \text{M}\right)  }$\textit{\ is good on M if
both }$\eta$\textit{\ }and $d\eta$ have Poincare growth\textit{. }Let
$\mathcal{E}$ be a vector bundle on $\mathcal{M}$(M) with a Hermitian metric
h. We will call h a good metric on $\overline{\mathcal{M}\left(
\text{M}\right)  }$ if the following holds: \textbf{i. }for all x$\in
\overline{\mathcal{M}\left(  \text{M}\right)  }\backslash\mathcal{M}\left(
\text{M}\right)  ,$ there exists sections $e_{1},...,e_{m}$ of $\mathcal{E}$
\ which forma basis of \textit{\ \ }$\mathcal{E}\left|  _{D^{r}\backslash
(D^{r}\cap\mathcal{D}_{\infty})}\right.  .$ \textbf{ii. }In a neighborhood
\textit{D}$^{r}$ of x in which $\overline{\mathcal{M}\left(  \text{M}\right)
}\backslash\mathcal{M}\left(  \text{M}\right)  $ is given by\textit{\ }%
$z_{1}\times.\times z_{k}=0.$ \textbf{iii. }The metric \textit{\ }%
h$_{i\overline{j}}=$h($e_{i},e_{j}$) has the following properties: \textbf{a.}
$\left|  h_{i\overline{j}}\right|  ,$ $\left(  \det\left(  h\right)  \right)
^{-1}\leq C\left(  \sum_{i=1}^{k}\log\left|  z_{i}\right|  \right)  ^{2m},$
\textit{for some }$C>0,$ $m\geq0.$ \textbf{b.}\ \textit{The 1-forms }$\left(
\left(  dh\right)  h^{-1}\right)  _{i\overline{j}}$ \textit{are good forms on
}$\overline{\mathcal{M}\left(  \text{M}\right)  }\cap D^{r}.$

It is easy to prove that there exists a unique extension $\overline
{\mathcal{E}}$ of $\mathcal{E}$ \ on $\overline{\mathcal{M}\left(
\text{M}\right)  },$ i.e. $\overline{\mathcal{E}}$ is defined locally
as\ holomorphic sections of $\mathcal{E}$ which have finite norm in h.

\begin{theorem}
\label{Mum100}Let ($\mathcal{E}$,h) be a vector bundle with a good metric on
$\mathcal{M}\left(  \text{M}\right)  $, then the Chern classes c$_{k}%
$($\mathcal{E}$,h) are good forms on $\overline{\mathcal{M}\left(
\text{M}\right)  }$ \textit{and the currents [c}$_{k}(\mathcal{E}$,h]
\textit{represent the cohomology classes }
\end{theorem}

\begin{center}
c$_{k}(\mathcal{E}$,h)$\in H^{2k}(\overline{\mathcal{M}\left(  \text{M}%
\right)  },\mathbb{Z)}.$
\end{center}

\subsubsection{PROOF\ OF\ THEOREM\ \ref{Mum100}}

For the proof see \cite{Mu}. $\blacksquare.$

\subsection{Applications of Mumford's Results to the Moduli of Odd Dimensional CY}

We are going to prove the following result:

\begin{theorem}
\label{Nik}Let $\pi:\mathcal{X\rightarrow}\mathcal{M}\left(  \text{M}\right)
$ be the flat family of non-singular CY manifolds. Let $\pi_{\ast}%
(\Omega_{\mathcal{X}\text{/}\mathcal{M}\text{(M)}}^{n,0})$ \textit{be equipped
with the metric h defined as follows:}
\end{theorem}

\begin{center}
$h(\omega_{\tau},\omega_{\tau})=(-1)^{\frac{n(n-1)}{2}}\left(  \sqrt
{-1}\right)  ^{n}\int_{\text{M}}\omega_{\tau}\wedge\overline{\omega_{\tau}}.$
\end{center}

\textit{Then h\ is a good metric.}

\subsubsection{PROOF\ OF\ THEOREM\ \ref{Nik}}

Let $\tau_{0}\in\mathcal{D}_{\infty}$ and let D be a one dimensional disk in
$\overline{\mathcal{M}\left(  \text{M}\right)  }$ which intersects
$\mathcal{D}_{\infty}$ on $\tau_{0}$ and D$^{\ast}=$D%
$\backslash$%
$\tau_{0}\in\mathcal{M}\left(  \text{M}\right)  .$ Over D%
$\backslash$%
$\tau_{0}$ we have a family $\mathcal{M}_{D^{\ast}}\mathcal{\rightarrow}$D%
$\backslash$%
$\tau_{0}$ of CY manifolds. We will assume that D is the unit disk and
$\tau_{0}$ is the origin of the disk. We know from the theory of Hodge
structures that if $\{\gamma_{1},...,\gamma_{b_{n}}\}$ is a basis of $H^{n}%
($M,$\mathbb{Z}$)/Tor then the functions:

\begin{center}
$\left(  ...,\int_{\gamma_{i}}\omega_{\tau},...\right)  $
\end{center}

for $0<\left|  \tau\right|  <1$ and $0<\arg(\tau)<2\pi$ are solutions of
differential equation with regular singularities. From the fact that the
solutions of any differential equation with regular singularities has a
logarithmic growth and

\begin{center}
$h(\omega_{\tau},\omega_{\tau})=$ $\left(  ...,\int_{\gamma_{i}}\omega_{\tau
},...\right)  \left(  <\gamma_{i},\gamma_{j}>\right)  \left(  ...,\int
_{\gamma_{i}}\overline{\omega_{\tau}},...\right)  ^{t}$
\end{center}

we deduce that

\begin{center}
$h(\omega_{\tau},\omega_{\tau})\leq C\left(  \sum_{i=1}^{k}\log\left|
\tau_{i}\right|  \right)  ^{2m}.$
\end{center}

From here we conclude that the form $\partial\left(  \log(h)\right)  $ has
also a logarithmic growth. Our theorem is proved. $\blacksquare.$

\subsection{Construction of the Analogue of the Dedekind $\eta$ Function for
Odd Dimensional CY manifolds}

\begin{theorem}
\label{Nik1}Let M be an odd dimensional CY manifold, then the C$^{\infty}$
section $\det(\overline{\partial})$ of the determinant line bundle
$\mathcal{L}$ constructed in Theorem \ref{CYsec} can be prolonged to a section
$\overline{\det(\overline{\partial})}$ of the line bundle $\overline
{\mathcal{L}}$ and $\overline{\det(\overline{\partial})}$ vanishes on the
discriminant locus $\mathcal{D}_{\infty}:=\overline{\mathcal{M}\text{(M)}%
}\backslash\mathcal{M}$(M).
\end{theorem}

\subsubsection{PROOF\ OF\ THEOREM\ \ref{Nik1}}

We know from Theorem \ref{CYsec} that the Quillen norm of the section
$\det(\overline{\partial})$ of the determinant line bundle $\mathcal{L}$ on
$\mathcal{M}$(M) is equal to the holomorphic Ray-Singer analytic torsion I(M),
i.e. $\left\|  \det(\overline{\partial})\right\|  _{Q}^{2}=$I(M). On the other
hand we proved in \cite{T1} that I(M)=$\left(  \det(\Delta_{0})\right)  ^{2}.$
This fact and Theorems \ref{var} and \ref{var1} imply that locally on
$\mathcal{M}$(M) the following formula is true:

\begin{center}
$dd^{c}\left(  \log\left(  \frac{I(M)}{\left\langle \omega_{\tau},\omega
_{\tau}\right\rangle }\right)  \right)  =dd^{c}\left(  \log\left(  \frac
{\det(\triangle_{0})}{\left\langle \omega_{\tau},\omega_{\tau}\right\rangle
}\right)  \right)  =0.$
\end{center}

From here we deduce in \cite{T1}\ that for each point $\tau\in\mathcal{M}(M)$
there exists an open set $\tau\in\mathcal{U}_{\tau}$ such that we have

\begin{center}
$\left\|  \det(\overline{\partial})\left|  _{\mathcal{U}_{\tau}}\right.
\right\|  _{Q}^{2}=\left\langle \omega_{\tau},\omega_{\tau}\right\rangle
\left|  f(\tau)\right|  ^{2},$
\end{center}

where $f$ is a holomorphic function in $\mathcal{U}_{\tau}$. Let us choose
$\mathcal{U}$ such that $\mathcal{U\cap}\overline{\mathcal{M}\text{(M)}}%
\neq\emptyset.$ Let $\tau_{0}\in\mathcal{D}_{\infty}\cap\mathcal{U}.$ We will
prove that we can continue $f$ locally around a point $\tau_{0}\in
\mathcal{D}_{\infty}=\overline{\mathcal{M}\text{(M)}}\backslash\mathcal{M}%
$(M). Indeed we proved in Theorem \ref{Nik} that $\left\langle \omega_{\tau
},\omega_{\tau}\right\rangle $ have a logarithmic growth around $\tau_{0}%
\in\mathcal{D}_{\infty}=\overline{\mathcal{M}\text{(M)}}\backslash\mathcal{M}%
$(M). Theorem \ref{LY} implies that $\left\|  \det(\overline{\partial
})\right\|  _{Q}^{2}=$I(M)$=\left(  \det(\triangle_{0})\right)  ^{2}\leq C$
and I(M)%
$\vert$%
$_{\mathcal{U}_{\tau}}=<\omega,\omega>|f|^{2}$ so we can conclude that
$\underset{\tau\rightarrow\tau_{0}}{\lim}f(\tau)=0$ since$\underset
{\tau\rightarrow\tau_{0}}{\lim}$ $\left\langle \omega_{\tau},\omega_{\tau
}\right\rangle =\infty$ and $\left\langle \omega_{\tau},\omega_{\tau
}\right\rangle $ has a logarithmitic growth. From here we deduce $,$ and $f$
can be continued around any point $\tau\in\mathcal{D}_{\infty}=\overline
{\mathcal{M}\text{(M)}}\backslash\mathcal{M}$(M). This fact and Theorem
\ref{Nik} implies that the canonical section $\det(\overline{\partial
})_{\mathcal{U}}$ can be prolonged to a section $\overline{\det(\overline
{\partial})}_{\mathcal{U}}$ of the line bundle $\overline{\mathcal{L}}$ and
$\overline{\det(\overline{\partial})}$ vanishes on the discriminant locus
$\mathcal{D}_{\infty}:=\overline{\mathcal{M}\text{(M)}}\backslash\mathcal{M}%
$(M). Theorem\ $\ref{Nik1}$ is proved. $\blacksquare.$

\begin{theorem}
\label{LY2}There exists a multi valued holomorphic section $\eta$ of the dual
of the determinant line bundle $\mathcal{L}$ over $\mathcal{M}(M)$ such that
the norm of $\eta$ with respect to the metric defined in Theorem \ref{Nik} is
equal to the Ray Singer Analytic torsion I(M)$=\det(\triangle_{0})^{2}.$
\end{theorem}

\subsubsection{PROOF\ OF\ THEOREM\ \ref{LY2}}

Theorem \ref{CYsec} implies that the line bundle $\mathcal{L}$ is a trivial
C$^{\infty}$ bundle. So the pullback $\pi^{\ast}(\mathcal{L})$ on the
universal cover $\widetilde{\mathcal{M}(M)}$ of $\mathcal{M}(M)$ will be a
trivial line bundle. Let $\widetilde{\sigma}$ be any non vanishing section of
$\pi^{\ast}(\mathcal{L}).$ Then since $\mathcal{L}=\widetilde{\mathcal{M}%
(M)}\times\mathbb{C}/\thicksim,$ where $(\tau,t)\thicksim(\tau_{1},t_{1})$
$\Leftrightarrow\tau_{1}=g\tau$ and $t_{1}=\chi(g)t,$ where $g\in\pi
_{1}(\mathcal{M}(M))$ and $\chi$ is a chararcter of $\pi_{1}(\mathcal{M}(M))$
that defines $\mathcal{L}.$ Theorem \ref{LY2} is proved. $\blacksquare.$

We know from the results in \cite{BGS1} and \ref{BGS2} that the determinant
line bundle $\mathcal{L}$ over $\mathcal{M}$(M) has a holomorphic structure.
In the next Theorem we will consider $\mathcal{L}$ as a holomorphic line
bundle defined over $\mathcal{M}$(M).

\begin{theorem}
\label{Nik2}There exists a positive integer $N$ such that $\mathcal{L}%
^{\otimes N}$ is a trivial complex analytic line bundle over $\mathcal{M}$(M).
\end{theorem}

\subsubsection{PROOF\ OF\ THEOREM\ \ref{Nik2}}

According to Theorem \ref{Tod1} $\mathcal{L}\approxeq R^{0}\pi_{\ast}%
(\Omega_{\mathcal{X}/\mathcal{M}(M)}^{2n+1,0}),$ where $\dim_{\mathbb{C}%
}M=2n+1.$ Therefore $\mathcal{L}$ is a subbundle of the flat vector bundle
$R^{2n+1}\pi_{\ast}\mathbb{C}$ where $\mathbb{C}$ is the locally constant
sheaf on $\mathcal{X},$ and $\pi:\mathcal{X\rightarrow M}(M)$ is the versal
family of CY manifolds over $\mathcal{M}(M)$. We know from Theorem \ref{Vie}
that $\mathcal{M}(M)=\mathcal{T}(M)/\Gamma,$ where $\mathcal{T}(M)$ is the
Teichm\"{u}ller space and $\Gamma$ is a subgroup of the mapping class group of
$M$ and according to \cite{Sul} $\Gamma$ is an arithmetic group.

If we lift the flat bundle $R^{n}\pi_{\ast}\mathbb{C}$ on $\mathcal{T}(M)$,
then $R^{2n+1}\pi_{\ast}\mathbb{C}$ will be a trivial bundle isomorphic to
$\mathcal{T}(M)\times H^{2n+1}(M_{0},\mathbb{C})$. Let us denote by

\begin{center}
$\sigma:\mathcal{T}(M)\rightarrow\mathcal{M}(M)=\mathcal{T}(M)/\Gamma$
\end{center}

the natural map. Clearly $\sigma^{\ast}(\mathcal{L})$ will be a flat complex
analytic subbundle of the trivial bundle $\mathcal{T}(M)\times H^{2n+1}%
(M_{0},\mathbb{C})$.

\begin{proposition}
\label{Nik21}Let N be a quasi-projective variety, $\mathcal{E}\approxeq
\mathbb{C}^{n}\times N$ be a trivial bundle and $\mathcal{L}$ be a flat line
bundle over N such that $\mathcal{L}^{\ast}\mathcal{\subset E},$ then
$\mathcal{L}$ is also trivial.
\end{proposition}

\textbf{Proof of Proposition }\ref{Nik21}: Let $\widetilde{N}$ be the
universal cover of $N.$ Clearly $\pi_{1}(N)$ acts without fixed points on
$\widetilde{N}.$ The pullback of $\mathcal{L}$ on $\widetilde{N}$ will be
denote by $\widetilde{\mathcal{L}}.$ $\widetilde{\mathcal{L}}$ will be a
trivial line bundle since $\mathcal{L}^{\ast}$ is a flat bundle over $N.$ Let
us denote by $\widetilde{\mathcal{E}}$ the pullback of $\mathcal{E}$ on
$\widetilde{\mathcal{N}}.$

$\pi_{1}(N)$ acts trivially on the trivial vector bundle $\widetilde
{\mathcal{E}}$ since $\mathcal{E}$ is a trivial vector bundle on $N.$ The
condition $\mathcal{L\subset E}$ implies that $\widetilde{\mathcal{L}}%
\subset\widetilde{\mathcal{E}}.$ $\pi_{1}(N)$ acts trivially on the trivial
line bundle $\widetilde{\mathcal{L}}$ since $\pi_{1}(N)$ acts trivially on the
trivial vector bundle $\widetilde{\mathcal{E}}$ and $\widetilde{\mathcal{L}%
}\mathcal{\subset}\widetilde{\mathcal{E}}.$ $\mathcal{L}$ is a trivial bundle
over $N$ since $\pi_{1}(N)$ acts trivially on the trivial line bundle
$\widetilde{\mathcal{L}}.$ Proposition \ref{Nik21} is proved. $\blacksquare.$

Proposition \ref{Nik21} implies we that $\sigma^{\ast}(\mathcal{L})$ will be a
trivial line bundle. The fact that $\sigma^{\ast}(\mathcal{L})\approxeq
\mathbb{C\times}\mathcal{T}(M)$ is a trivial bundle implies that
$\mathcal{L}\approxeq\mathbb{C\times}\mathcal{T}(M)/\Gamma,$ where $\Gamma$
acts in a natural way on the Teichm\"{u}ller space and it acts by a character
$\chi\in Hom(\Gamma,\mathbb{C}_{1}^{\ast})\approxeq Hom(\Gamma/[\Gamma
,\Gamma],\mathbb{C}_{1}^{\ast})$ on $\mathbb{C}.$ According to a Theorem of
Kazhdan $\Gamma/[\Gamma,\Gamma]$ is a finite group since $\Gamma$ is an
arithmetic group of rank $\geq2$ according to \cite{Sul} and \cite{Bour}. From
here we deduce that $\mathcal{L}^{N}$ will be a trivial bundle on
$\mathcal{M}(M)$, where $N=\#\Gamma/[\Gamma,\Gamma].$ Theorem \ref{Nik2} is
proved. $\blacksquare.$

\begin{corollary}
\label{Nik22} There exists a holomorphic section $\eta^{N}$ of the trivial
bundle $\left(  \mathcal{L}^{\ast}\right)  ^{N}$ such that it can be prolonged
to a holomorphic section $\overline{\eta^{N}}$ \textit{of (}$\overline
{\mathcal{L}^{\ast}})^{N}$ \textit{whose zero set is supported by
}$\mathcal{D}_{\infty}$ and \textit{the Quillen norm }$\left\|  \eta
^{N}\right\|  _{Q}^{2}=\det(\Delta_{0})^{2N}.$
\end{corollary}

\textbf{Proof of Corollary }\ref{Nik22}\textbf{: }As we pointed out we can
prolonged the dual of the determinant line bundle $\mathcal{L}$ from
$\mathcal{M}$(M) to a holomorphic line bundle$\mathcal{\ }\overline
{\mathcal{L}^{\ast}}$ over $\overline{\mathcal{M}(M)}.$ In Theorem \ref{CYsec}
we constructed a non vanishing $C^{\infty}$ section $\det(\overline{\partial
})$ of $\mathcal{L}^{\ast}$ over $\mathcal{M}$(M). We know that the norm of
$\det(\overline{\partial})$ with respect of the metric on $\mathcal{L}$
defined in Theorem \ref{Nik} is exactly equal to the Ray Singer Analytic
Torsion I(M)$=\det(\Delta_{0})^{2}$. From Theorem \ref{Nik1} we know that we
can prolong the section $\det(\overline{\partial})$ to a section
$\overline{\det(\overline{\partial})}\ $of the line bundle $\overline
{\mathcal{L}^{\ast}}$ over $\overline{\mathcal{M}(M)}$ and the support of the
zero set of $\overline{\det(\overline{\partial})}$ is exactly $\mathcal{D}%
_{\infty}.$ From here we deduce that the Poincare dual homology class of the
Chern class of $\overline{\mathcal{L}\text{ }}$ is an effective divisor whose
support is the same as that of $\mathcal{D}_{\infty}.$ Combining this fact
together with Theorem \ref{Nik2} we conclude that there exists a holomorphic
section $\overline{\mathcal{\eta}^{N}\text{ }}$ of the line bundle $\left(
\overline{\mathcal{L}^{\ast}}\right)  ^{N}$ over $\overline{\mathcal{M}(M)}$
whose zero set is supported by $\mathcal{D}_{\infty}$ and the multiplicities
of the components of the irreducible divisors of $\left(  \overline
{\mathcal{\eta}}\right)  ^{N}$ are the same as of $\left(  \overline
{\det(\overline{\partial})}\right)  ^{N}.$ The fact that $\eta^{N\text{ }}$ is
defined up to a constant and Corollary \ref{LY2} we conclude that after
normalizing $\eta^{N}$ its Quillen norm $\left\|  \eta^{N}\right\|  _{Q}^{2}$
will be $\det(\Delta_{0})^{2N}.$ Corollary \ref{Nik22}\textbf{ }is proved.
$\blacksquare$

\section{Some Problems}

Let us define Sh(C,E,Z) to be the set of all families $\pi:Y\rightarrow C$ of
fixed type CY manifolds Z over a fix complete algebraic curve C with a fixed
set of points over which the fibres are singular, up to isomorphisms.

\begin{problem}
\label{pr0}Is the set Sh(C,E,Z) finite?
\end{problem}

The results of this paper combined with the results from \cite{JT97} imply
that the set Sh(C,E,Z) is discrete. In order to prove that it is finite one
need to find a uniform bound on the volume of the image of the curve in the
moduli space $\mathcal{M(}M)$ of CY manifolds. We found a bound of the images
of fix C and fix set of points E on C over we which the fibres are singular by
using Gauss-Bonnet theorem and the fact that Weil-Petersson metric is complete
on the moduli space of pseudo polarized algebraic K3 surfaces. This method
does not work for CY threefolds, since Weil-Petersson metric is not complete.

\begin{problem}
\label{pr}Can one find a product formulas for $\eta$ around points of maximal
degenerations, which means that around \ such points the monodromy operator
has index of unipotency $n+1?$
\end{problem}

For more precise discussion of Problem \ref{pr} see \cite{T2}. Problem
\ref{pr} is closely related to the paper \cite{BCOV} and more precisely to the
counting problem of elliptic curves on CY threefold.

\begin{problem}
\label{pr1}Is is true that any CY manifold can be deformed to an algebraic
manifold with one conic singularity?
\end{problem}

Problem \ref{pr1} is related to the following problem: Let M be an algebraic
variety embedded in $\mathbb{P}^{N}.$ Suppose that the component of the
Hilbert scheme $\mathcal{H}_{M/\mathbb{P}^{N}}^{\prime}$ that contains an open
non singular quasi-projective subscheme $\mathcal{H}_{M/\mathbb{P}^{N}}. $ Let
$\mathcal{D}_{M}$ be the set of the points in $\mathcal{H}_{M/\mathbb{P}^{N}}$
that corresponds to singular varieties. It is not difficult to prove that
$\mathcal{D}_{M}$ is a closed subvariety in $\mathcal{H}_{M/\mathbb{P}^{N}}.$
Suppose that $\mathcal{D}_{M}$ is a divisor in $\mathcal{H}_{M/\mathbb{P}^{N}}.$

\begin{problem}
\label{pr2}Is it true that the generic point of $\mathcal{D}_{M}$ corresponds
to a projective manifold with only one conic singularity?
\end{problem}

The results of this paper suggest that one can expect that the Hilbert scheme
$\mathcal{H}_{M/\mathbb{P}^{N}}^{\prime}$ of odd dimensional CY manifold M
contains a non singular open quasi-projective variety $\mathcal{H}%
_{M/\mathbb{P}^{N}}$ such that the discriminant locus $\mathcal{D}_{M}$ is a divisor.

Problem \ref{pr2} is closely related to the Miles Ried's conjecture that the
moduli of all CY threefolds is connected.


\begin{thebibliography}{99}
\bibitem{ABKS}D. Abramovich, J. -F. Burnol, J. Kramer and C. Soul\'{e},
''\textit{Lectures on Arakelov Geometry''}, Cambridge Studies In Advanced
Mathematics Volume \textbf{33}, Cambridge University Press, 1992.

\bibitem{BCOV}M. Bershadsky, S. Cecotti, H. Ooguri and C. Vafa,
\textit{''Kodaira-Spencer Theory of Gravity and Exact Results for Quantum
String Amplitude''}, Comm. Math. Phys. \textbf{165} (1994), 311-428.

\bibitem{BF}J.-M. Bismut and D. Freed, ''\textit{The Analysis of Elliptic
Families I. Metrics and Connections on Determinant Line Bundles'',} Comm.
Math. Phys. \textbf{106,} 159-176(1986).

\bibitem{Bi}E. Bishop, \textit{''Conditions for the Analyticity of Certain
Sets.'' }Mich. Math. J. \textbf{11}(1964) 289-304.

\bibitem{BGS1}J.-M. Bismut, H. Gillet and C. Soul\`{e}, ''\textit{Analytic
Torsion and Holomorphic Determinant Bundles I. Bott-Chern Forms and Analytic
Torsion'', }Comm. Math. Phys.\textbf{115,} 49-78(1988).

\bibitem{BGS2}J.-M. Bismut, H. Gillet and C. Soul\`{e}, ''\textit{Analytic
Torsion and Holomorphic Determinant Bundles II. Direct Images and Bott-Chern
Forms '', }Comm. Math. Phys.\textbf{115,} 79-126(1988).

\bibitem{BGV}N. Berline, E. Getzler and \ M. Vergne, \textit{''Heat Kernels
and Dirac Operators'', }Springer Verlag, 1991.

\bibitem{Bour}A. A. Kirillov and Cl. Delaroche, ''\textit{Sur les Relations
entre l'espace dual d'un group et la struture de ses sous-groupes ferm\'{e}s''
(d'apr\`{e}s D. A. Kazhdan)}, Bourbaki No. 343, 1967/68.

\bibitem{D}S. Donaldson, ''\textit{Infinite Determinants, Stable Bundles and
Curvature'', }Duke Mathematical Journal, vol. \textbf{54, Number 1}(1987), 231-248.

\bibitem{DK}S. Donaldson and P. Kronheimer, ''\textit{The Geometry of Four
Manifolds'',} Oxford Mat. Monog., Oxford Science Publications, Oxford
University Press, New York 1990.

\bibitem{F}D. Freed, ''\textit{On Determinant Line Bundles'',} Mathematical
Aspects of String Theory, ed. S.-T. Yau, World Scientific, Singapore, New
Jersey, Hong Kong, p. 189-238(1998).

\bibitem{Gil}P. Gilkey, ''\textit{Invariance Theory, The Heat Equation, And
the Atiyah-Singer Index Theorem'',} Mathematics Lecture Series vol.
\textbf{11}, Publish or Parish, Inc. Wilmington, Delaware (USA) 1984.

\bibitem{Gr}Ph. Griffiths, \textit{''Periods of Integrals on Algebraic
Manifolds'' I and II,} Amer. Jour. Math. \textbf{90} (1968), 568-626 and 805-865.

\bibitem{JT95}J. Jorgenson and A. Todorov, ''\textit{A Conjectural Analogue of
Dedekind Eta Function for K3 Surfaces'', }\ Math. Research Lett.
\textbf{2}(1995) 359-360.

\bibitem{JT96}J. Jorgenson and A. Todorov, ''\textit{Analytic Discriminant for
Manifolds with Zero Canonical Class'', }Manifolds and Geometry, ed. P. de
Bartolomeis, F. Tircerri and E. Vesantini, Symposia Mathematica \textbf{36,
}(1996) 223-260.

\bibitem{JT098}J. Jorgenson and A. Todorov, in preparation.

\bibitem{JT97}J. Jorgenson and A. Todorov, ''\textit{Ample Divisors,
Automorphic Forms and Shafarevich's Conjecture''}, preprint 1999.

\bibitem{JT98}J. Jorgenson and A. Todorov, ''\textit{Analytic Discriminant for
Polarized Algebraic K3 Surfaces'',} Mirror Symmetry III, ed. S-T. Yau and
Phong, AMS, p. 211-261.

\bibitem{KM}K. Kodaira and Morrow, ''\textit{Complex Manifolds''.}

\bibitem{MM}T. Matsusaka and D. Mumford, ''\textit{Two Fundamental Theorems on
Deformations of Polarized Varieties'',} Amer. J. Math. \textbf{86}(1964).

\bibitem{Mu}D. Mumford, ''\textit{Hirzebruch's Proportionality Principle in
the Non-Compact Case'',} Inv. Math. \textbf{42}(1977), 239-272.

\bibitem{Roe}J.Roe, ''\textit{Elliptic Operators, Topology and Asymptotic
Methods''} Pitman Research Notes in Mathematics Series \textbf{179,} Longman
Scientific \& Technical, \ 1988.

\bibitem{RS}D. Ray and I. Singer, ''\textit{Analytic Torsion for Complex
Manifolds'',} Ann. Math. \textbf{98 }(1973) 154-177.

\bibitem{Sul}D. Sullivan, ''\textit{Infinitesimal Computations in Topology'',
Publ. Math. IHES, No \textbf{47 (}1977), 269-331.}

\bibitem{Ti}G. Tian, ''\textit{Smoothness of the Universal Deformation Space
of Calabi-Yau Manifolds and its Petersson-Weil Metric'',} Math. Aspects of
String Theory, ed. S.-T.Yau, World Scientific (1998), 629-346.

\bibitem{To89}A. Todorov, \textit{''The Weil-Petersson Geometry of \ Moduli
Spaces of SU(n}$\geq3$\textit{) (Calabi-Yau Manifolds) I'',} Comm. Math. Phys.
\textbf{126} (1989), 325-346.

\bibitem{TO99}A. Todorov, \textit{''Witten's Geometric Quantization of Moduli
of CY manifolds'',} preprint 1999.

\bibitem{T1}A. Todorov, \textit{''Ray Singer Analytic Torsion of CY\ Manifolds
I''}. Preprint.

\bibitem{T2}A. Todorov, ''\textit{Weil-Petersson Geometry of CY\ Manifolds
II',} preprint.

\bibitem{W}E. Viehweg, \textit{''Quasi-Projective Moduli for Polarized
Manifolds'', }\ Ergebnisse der Mathematik und iher Grenzgebiete 3. Folge, Band
\textbf{30}, Springer-Ver;ag, 1991.

\bibitem{Y}Yoshikawa, ''\textit{Generalized Enriques Surfaces and Analytic
Torsion'', }preprint (1998).

\bibitem{Yau}S. T. Yau, ''\textit{On the Ricci Curvature of Compact K\"{a}hler
Manifolds and Complex Monge-Amper Equation I''}, Comm. Pure and App. Math.
\textbf{31,} 339-411(1979).
\end{thebibliography}
\end{document}